\documentclass[12pt]{article}
\usepackage{amsmath, amsthm, amssymb}

\theoremstyle{plain}
\newtheorem{thm}{Theorem}
\newtheorem*{MainLemma}{The Main Lemma}
\newtheorem*{KostochkasLemma}{Kostochka's Lemma}
\newtheorem*{CapRizz}{Caprara and Rizzi}
\newtheorem*{FractionalTheorem}{Fractional Version}
\newtheorem*{KingVettaReed}{King, Reed and Vetta}
\newtheorem*{ReedConjecture}{Reed's Conjecture}
\newtheorem*{TheoremD}{Theorem D}

\newtheorem{lem}[thm]{Lemma}

\theoremstyle{definition}

\newtheorem*{CliqueGraph}{Clique Graph}
\newtheorem*{Satisfaction}{Satisfaction}
\theoremstyle{remark}

\newtheorem*{observation}{Observation}

\title{On hitting all maximum cliques with an independent set}
\author{Landon Rabern\\
\small \texttt{landon.rabern@gmail.com}}
\begin{document}
\maketitle

\begin{abstract}
We prove that every graph $G$ for which $\omega(G) \geq \frac{3}{4}(\Delta(G) + 1)$, has an independent set $I$ such that $\omega(G - I) < \omega(G)$.  It follows that a minimum counterexample $G$ to Reed's conjecture satisfies $\omega(G) < \frac{3}{4}(\Delta(G) + 1)$ and hence also $\chi(G) > \left\lceil \frac{7}{6}\omega(G) \right\rceil$.  We also prove that if for every induced subgraph $H$ of $G$ we have $\chi(H) \leq \max\left\{\left\lceil \frac{7}{6}\omega(H) \right\rceil, \left\lceil \frac{\omega(H) + \Delta(H) + 1}{2}\right\rceil\right\}$, then we also have $\chi(G) \leq \left\lceil \frac{\omega(G) + \Delta(G) + 1}{2}\right\rceil$.  This gives a generic proof of the upper bound for line graphs of multigraphs proved by King et al.  
\end{abstract}

\section{Introduction}
\begin{Satisfaction}
Given a relation $R$ between graph functions we say that a graph $G$ \emph{satisfies} $R$ if plugging $G$ into each function gives a true statement (e.g. $G$ satisfies $\chi \geq \Delta$ means $\chi(G) \geq \Delta(G)$).
\end{Satisfaction}

We prove the following general lemma and apply it to Reed's conjecture.

\begin{MainLemma}
If $G$ is a graph satisfying $\omega \geq \frac{3}{4}(\Delta + 1)$, then $G$ has an independent set $I$ such that $\omega(G - I) < \omega(G)$.
\end{MainLemma}

In \cite{Reed}, Reed conjectured the following upper bound on the chromatic number.

\begin{ReedConjecture}
Every graph satisfies $\chi \leq \left\lceil \frac{\omega + \Delta + 1}{2}\right\rceil$.
\end{ReedConjecture}

\begin{observation}
If we could always find an independent set whose removal decreased both $\omega$ and $\Delta$, then the conjecture would follow by simple induction.  Expanding the independent set given by The Main Lemma to a maximal one shows that this sort of argument goes through when $\omega \geq \frac{3}{4}(\Delta + 1)$.
Thus a minimum counterexample to Reed's conjecture satisfies $\omega < \frac{3}{4}(\Delta + 1)$ and hence also $\chi > \left\lceil \frac{7}{6}\omega \right\rceil$.\\
\end{observation}

Reed's upper bound was proved for line graphs of multigraphs by King, Reed and Vetta in \cite{LineGraphs}, for quasi-line graphs by King and Reed in \cite{QuasiLineGraphs}, and recently King and Reed proved it for all claw-free graphs (see King's thesis \cite{KingThesis}).  The line graphs of multigraphs result follows from the following theorem.\newline

\begin{TheoremD}
If every induced subgraph of a graph $G$ satisfies $\chi \leq \max\left\{\left\lceil \frac{7}{6}\omega \right\rceil, \left\lceil \frac{\omega + \Delta + 1}{2}\right\rceil\right\}$, then $G$ satisfies $\chi \leq \left\lceil \frac{\omega + \Delta + 1}{2}\right\rceil$.
\end{TheoremD}

Reed's upper bound for line graphs of multigraphs follows immediately from Theorem D, a bound of Caprara and Rizzi (see \cite{CapraraAndRizzi}) and the bound of Molloy and Reed on fractional colorings (see \cite{MolloyAndReed}).

\begin{CapRizz}
Let $H$ be a multigraph and $G = L(H)$.  Then
\[\chi(G) \leq \max\left\{\left\lfloor 1.1\Delta(H) + 0.7 \right\rfloor, \left\lceil \chi^*(G) \right\rceil\right\}\]
\end{CapRizz}

\begin{FractionalTheorem}
Every graph satisfies $\chi^* \leq \frac{\omega + \Delta + 1}{2}$.
\end{FractionalTheorem}

Since induced subgraphs of line graphs are line graphs, $\Delta(H) \leq \omega(G)$ and $\left\lfloor 1.1 \omega(G) + 0.7 \right\rfloor \leq \left\lceil \frac{7}{6}\omega(G) \right\rceil$ we may apply Theorem D and prove the following.

\begin{KingVettaReed}
If $G$ is the line graph of a multigraph, then $G$ satisfies $\chi \leq \left\lceil \frac{\omega + \Delta + 1}{2}\right\rceil$.
\end{KingVettaReed}

\section{Proof of The Main Lemma}
We need three lemmas.  The first is due to Hajnal (see \cite{Hajnal}).

\begin{lem}\label{HajnalLemma}
Let $G$ be a graph and $\mathcal{Q}$ a collection of maximum cliques in $G$. Then
\[\left | \bigcap \mathcal{Q}\right | \geq 2\omega(G) - \left | \bigcup \mathcal{Q}\right |.\]
\end{lem}
\begin{proof}
Assume (to reach a contradiction) that the lemma is false and let $\mathcal{Q}$ be a counterexample with $|\mathcal{Q}|$ minimal.  Put $r = |\mathcal{Q}|$ and $\mathcal{Q} = \{Q_1, ..., Q_r\}$.  Consider the set $\displaystyle W = (Q_1 \cap \bigcup_{i=2}^r Q_i) \cup \bigcap_{i=2}^r Q_i$.  Plainly, $W$ is a clique.  Thus
\begin{align*}
\omega(G) &\geq |W| \\
&= \left |(Q_1 \cap \bigcup_{i=2}^r Q_i) \cup \bigcap_{i=2}^r Q_i \right | \\
&= \left |Q_1 \cap \bigcup_{i=2}^r Q_i \right | + \left |\bigcap_{i=2}^r Q_i \right | - \left |\bigcap_{i=1}^r Q_i \cap \bigcup_{i=2}^r Q_i \right | \\
&= \left |Q_1 \right | + \left |\bigcup_{i=2}^r Q_i \right | - \left |\bigcup_{i=1}^r Q_i \right| + \left |\bigcap_{i=2}^r Q_i \right | - \left |\bigcap_{i=1}^r Q_i \right | \\
&= \omega(G) + \left |\bigcup_{i=2}^r Q_i \right | + \left |\bigcap_{i=2}^r Q_i \right | - \left |\bigcup_{i=1}^r Q_i \right| - \left |\bigcap_{i=1}^r Q_i \right | \\
&\geq \omega(G) + 2\omega(G) - \left |\bigcup_{i=1}^r Q_i \right| - \left |\bigcap_{i=1}^r Q_i \right |. \\
\end{align*}

Thus $\displaystyle \left |\bigcap_{i=1}^r Q_i \right | \geq 2\omega(G) - \left |\bigcup_{i=1}^r Q_i \right|$ giving a contradiction.
\end{proof}

The second lemma we need is an improvement of Hajnal's result for graphs satisfying $\omega > \frac{2}{3}(\Delta + 1)$ due to Kostochka (see \cite{Kostochka}).  We give a new (simpler) proof of this result.

\begin{CliqueGraph}
Let $G$ be a graph. For a collection of cliques $\mathcal{Q}$ in $G$, let $X_{\mathcal{Q}}$ be the intersection graph of $\mathcal{Q}$.  That is, the vertex set of $X_{\mathcal{Q}}$ is $\mathcal{Q}$ and there is an edge between $Q_1 \neq Q_2 \in \mathcal{Q}$ if and only if $Q_1$ and $Q_2$ intersect.
\end{CliqueGraph}

\begin{lem}\label{TwoThirdsLemma}
Let $G$ be a graph satisfying $\omega > \frac{2}{3}(\Delta + 1)$.   If $\mathcal{Q}$ is a collection of maximum cliques in $G$ such that $X_{\mathcal{Q}}$ is connected, then $X_{\mathcal{Q}}$ is complete.
\end{lem}
\begin{proof}
Let $Q_1, Q_2, Q_3 \in \mathcal{Q}$ be distinct and assume that $Q_1 \cap Q_2 \neq \emptyset$ and $Q_2 \cap Q_3 \neq \emptyset$.  Then $|Q_1 \cap Q_2| = |Q_1| + |Q_2| - |Q_1 \cup Q_2| \geq 2\omega(G) - (\Delta(G) + 1)$.  Hence

\begin{align*}
|Q_1 \cap Q_3| &\geq |Q_1 \cap Q_2 \cap Q_3| \\
&\geq |Q_1 \cap Q_2| - (|Q_2| - |Q_2 \cap Q_3|) \\
&\geq 2\omega(G) - (\Delta(G) + 1) - (\omega(G) - (2\omega(G) - (\Delta(G) + 1))) \\
&= 3\omega(G) - 2(\Delta(G) + 1) > 0.\\
\end{align*}

Thus $Q_1 \cap Q_3 \neq \emptyset$ showing that $X_{\mathcal{Q}}$ is transitive.  The lemma follows since a transitive connected graph is complete.
\end{proof}

\begin{KostochkasLemma}
Let $G$ be a graph satisfying $\omega > \frac{2}{3}(\Delta + 1)$.   If $\mathcal{Q}$ is a collection of maximum cliques in $G$ such that $X_{\mathcal{Q}}$ is connected, then $\cap \mathcal{Q} \neq \emptyset$.
\end{KostochkasLemma}
\begin{proof}
Assume not and let $\mathcal{Q} = \{Q_1, ..., Q_r\}$ be a bad collection of maximum cliques with $r$ minimal.  Then $r \geq 3$ and $X_{\mathcal{Q}}$ is complete by Lemma \ref{TwoThirdsLemma}.  Put $\mathcal{Z} = \mathcal{Q} - \{Q_1\}$.  Then $X_{\mathcal{Z}}$ is connected and hence by minimality of $r$, we see that $\cap \mathcal{Z} \neq \emptyset$.  In particular $|\cup \mathcal{Z}| \leq \Delta(G) + 1$.  Thus

\[|\cup Q| \leq |Q_1 - Q_2| + |\cup \mathcal{Z}| \leq 2(\Delta(G) + 1) - \omega(G) < 2\omega(G)\]

But then Lemma \ref{HajnalLemma} gives a contradiction.
\end{proof}

The third lemma we need is a result of Haxell (see \cite{Haxell}) on independent transversals.

\begin{lem}\label{HaxellLemma}
Let $k$ be a positive integer, let $H$ be a graph of maximum degree at most $k$,
and let $V(H) = V_1 \cup \cdots \cup V_n$ be a partition of the vertex set of $H$. Suppose that $|V_i| \geq 2k$ for each $i$. Then $H$ has an independent set $\{v_1, \ldots, v_n\}$ where $v_i \in V_i$ for each $i$.
\end{lem}

\begin{proof}[Proof of The Main Lemma]
Let $G$ be a graph satisfying $\omega \geq \frac{3}{4}(\Delta + 1)$. Let $\mathcal{Q}$ be the collection of all maximum cliques in $G$ and let $\{\mathcal{C}_1, \ldots, \mathcal{C}_r\}$ be the vertex sets of the components of $X_{\mathcal{Q}}$. For each $i$, put $F_i = \bigcap \mathcal{C}_i$. Then, by Lemma \ref{HajnalLemma} and Kostochka's Lemma, we have $|F_i| \geq 2\omega(G) - (\Delta(G) + 1)$.  Since every vertex $v \in F_i$ is in a maximum clique in $\bigcup \mathcal{C}_i$, $v$ hits at most $\Delta(G) + 1 - \omega(G)$ vertices outside of $\bigcup \mathcal{C}_i$.\newline

Let $H$ be the graph with $V(H) = \displaystyle \bigcup_{i} F_i$ and an edge between $v, w \in V(H)$ if and only if $vw \in E(G)$ and $v$ and $w$ are in different $F_i$.  Then, by the above, $\Delta(H) \leq \Delta(G) + 1 - \omega(G)$.\newline

Consider the partition $\{F_i\}_{i}$ of $V(H)$.  We have 
\begin{align*}
|F_i| &\geq 2\omega(G) - (\Delta(G) + 1) \\
&\geq 2 \frac{3}{4}(\Delta(G) + 1) - (\Delta(G) + 1) \\
&= \frac{1}{2}(\Delta(G) + 1) \\
&= 2(\Delta(G) + 1 - \frac{3}{4}(\Delta(G) + 1)) \\
&\geq 2(\Delta(G) + 1 - \omega(G)) \\
&\geq 2\Delta(H). \\
\end{align*}

Thus, by Lemma \ref{HaxellLemma}, $H$ has an independent set $I = \{v_1, \ldots, v_n\}$ where $v_i \in F_i$ for each $i$.  Since $F_i$ was contained in each element of $\mathcal{C}_i$ we have $\omega(G - I) < \omega(G)$.
\end{proof}
\section{Proof of Theorem D}
Theorem D is an easy consequence of The Main Lemma.
\begin{proof}[Proof of Theorem D]
Assume (to reach a contradiction) that the theorem is false and let $G$ be a counterexample with the minimum number of vertices.  First assume that $\omega(G) \geq \frac{3}{4}(\Delta(G) + 1)$.  Then by The Main Lemma we have an independent set $I$ with $\omega(G - I) < \omega(G)$.  Plainly, we may assume that $I$ is maximal (and hence $\Delta(G - I) < \Delta(G)$. Put $H = G - I$.  Then, by minimality of $G$, we have
\begin{align*}
\chi(G) &\leq 1 + \chi(H) \\
&\leq 1 + \left\lceil \frac{\omega(H) + \Delta(H) + 1}{2}\right\rceil \\
&\leq 1 + \left\lceil \frac{\omega(G) - 1 + \Delta(G) - 1 + 1}{2}\right\rceil \\
&\leq \left\lceil \frac{\omega(G) + \Delta(G) + 1}{2}\right\rceil. \\
\end{align*}

This is a contradiction, hence we must have $\omega(G) < \frac{3}{4}(\Delta(G) + 1)$.  But then
\begin{align*}
\left\lceil \frac{7}{6}\omega(G) \right\rceil &\geq \chi(G) \\
&> \left\lceil \frac{\omega(G) + \Delta(G) + 1}{2}\right\rceil \\
&\geq \left\lceil \frac{\omega(G) + \frac{4}{3}\omega(G)}{2}\right\rceil \\
&= \left\lceil \frac{7}{6}\omega(G) \right\rceil. \\
\end{align*}

This final contradiction completes the proof.
\end{proof}


\begin{thebibliography}{1}

\bibitem{CapraraAndRizzi}
A. Caprara and R. Rizzi.
\newblock Improving a family of approximation algorithms to edge color multigraphs.
\newblock{\em Information Processing Letters,} \textbf{68}, 1998, \mbox{11 - 15.}

\bibitem{Hajnal}
A. Hajnal.
\newblock A theorem on k-saturated graphs. 
\newblock{\em Can. J. Math.,} \textbf{10(4)}, 1965, \mbox{720-724}.

\bibitem{Haxell}
P. E. Haxell.
\newblock A Note on Vertex List Colouring. 
\newblock{\em Combinatorics, Probability and Computing,} \textbf{10(4)}, 2001, \mbox{345-347}.

\bibitem{LineGraphs}
A. King, B. Reed, and A. Vetta. 
\newblock An upper bound for the chromatic number of line graphs. 
\newblock{\em European Journal of Combinatorics,} \textbf{28(8)}, 2007, \mbox{2182-2187}.

\bibitem{QuasiLineGraphs}
A. King, B. Reed.
\newblock Bounding $\chi$ in terms of $\omega$ and $\Delta$ for quasi-line graphs.
\newblock{\em Journal of Graph Theory,} To Appear.

\bibitem{KingThesis}
A. King.
\newblock Claw-free graphs and two conjectures on $\omega$, $\Delta$, and $\chi$.
\newblock{Thesis.}

\bibitem{Kostochka}
A. Kostochka. 
\newblock Degree, clique number, and chromatic number. 
\newblock{\em Metody Diskret. Anal.,} \textbf{35}, 1980, \mbox{45-70}. [In Russian]

\bibitem{MolloyAndReed}
M. Molloy and B. Reed.
\newblock Graph Coloring and the Probabilistic Method.
\newblock{\em Springer-Verlag, Berlin,} 2000.

\bibitem{Reed}
B. Reed.
\newblock $\omega$, $\Delta$, and $\chi$.
\newblock{\em Journal of Graph Theory,} \textbf{27}, 1998, \mbox{177-212}.

\end{thebibliography}
\end{document}